\documentclass[9pt,shortpaper,twoside,web]{ieeecolor}
\usepackage{generic}
\usepackage{cite}
\usepackage{amsmath,amssymb,amsfonts}
\usepackage{algorithmic}
\usepackage{graphicx}
\usepackage{pgfplots}
\usetikzlibrary{positioning}
\pgfplotsset{compat=1.14} 
\tikzstyle{vertex} = [fill,shape=circle,node distance=30pt]
\tikzstyle{edge} = [fill,opacity=.6,fill opacity=.5,line cap=round, line join=round, line width=10pt]
\tikzstyle{elabel} =  [fill,shape=circle,node distance=30pt,fill opacity=.9]
\definecolor{mygray}{gray}{0.95}
\pgfdeclarelayer{background}
\pgfsetlayers{background,main}
\usepackage[most]{tcolorbox}
\definecolor{mypurple}{rgb}{0.59, 0.44, 0.84}
\usepackage[mathscr]{euscript}

\usepackage{textcomp}

\begin{document}
\title{On the Stability of Multilinear Dynamical Systems}
\author{Can Chen
\thanks{This work is supported in part by Rackham Graduate Fellowship, in part by Air Force Office of Scientific Research Award No: FA9550-18-1-0028, and in part by  National Science Foundation Grant DMS 1613819. }
\thanks{Can Chen is with the Department of Mathematics, University of Michigan, Ann Arbor, MI 48109 USA and with Channing Division of Network Medicine, Brigham and Women's Hospital and Harvard Medical School, Boston, MA 02115 USA (canc@umich.edu)}}

\maketitle

\begin{abstract}
This paper investigates the stability properties of discrete-time multilinear dynamical systems via tensor spectral theory. In particular, if the dynamic tensor of a multilinear dynamical system is orthogonally decomposable (odeco), we can construct its explicit solution by exploiting tensor Z-eigenvalues and Z-eigenvectors. Based on the form of the explicit solution, we illustrate that the Z-eigenvalues of the dynamic tensor  play a significant role in the stability analysis, offering necessary and sufficient conditions. In addition, by utilizing the upper bounds of Z-spectral radii, we are able to determine the asymptotic stability of the multilinear dynamical system efficiently. Furthermore, we extend the stability results to the multilinear dynamical systems with non-odeco dynamic tensors by exploiting tensor singular values.  We demonstrate our results via numerical examples.
\end{abstract}

\begin{IEEEkeywords}
Multilinear dynamical systems, regions of attraction, stability, tensors, Z-eigenvalues, singular values. 
\end{IEEEkeywords}

\section{Introduction}
The role of tensor algebra has been explored for modeling and simulation of linear and nonlinear dynamics \cite{doi:10.1137/1.9781611975758.18,chen2021controllability, chen2021multilinear,gelss2017tensor,klus2018tensor,kruppa2017comparison}. The key idea is to tensorize the vector-based dynamical system representation into an equivalent tensor representation and to exploit tensor algebra. Tensor decomposition techniques such that CANDECOMP/PARAFAC decomposition, higher-order singular value decomposition, and tensor train decomposition, are applied for reducing memory usages and enabling efficient computations in multilinear dynamical systems theory \cite{chen_2020, chen2021multilinear,gelss2017tensor}. Gel{\ss} \cite{gelss2017tensor} applied tensor decompositions for computing numerical solutions of master equations associated with Markov processes on extremely large state spaces, developing tensor representations for operators based on nearest-neighbor interactions, construction of pseudoinverses for dimensionality reduction methods, and the approximation of transfer operators of dynamical systems. In addition, Chen \textit{et al.} \cite{doi:10.1137/1.9781611975758.18,chen2021multilinear} developed the tensor algebraic conditions for stability, reachability, and observability for input/output multilinear time-invariant systems, and expressed them in terms of more standard notions of tensor ranks/decompositions to facilitate efficient computations.

Many complex systems such as those arising in biology and engineering can be studied using a network prospective \cite{scale_free,chen2021controllability,dynamic_network,MONNIG2018347,Ranshous:2015:ADD:3160212.3160218,RIED20171,small_world}. Most real world data representations are multidimensional, and using graph models to characterize them may cause a loss of higher-order information \cite{9119161,chen2021controllability,7761624}. Recently, a new tensor-based continuous-time multilinear dynamical system representation with linear control inputs (different from the ones proposed in \cite{doi:10.1137/1.9781611975758.18,chen2021multilinear} which can be unfolded to linear dynamical systems via tensor unfolding, an operation that transforms a tensor into a matrix) was proposed by Chen \textit{et al.} \cite{chen2021controllability} for characterizing the multidimensional state dynamics of  hypergraphs, a generalization of graphs in which edges can contain more than one nodes. The authors derived a Kalman-rank-like condition to determine the minimum number of control nodes needed to achieve controllability/accessibility of  hypergraphs. They also proposed minimum number of control nodes as a measure of hypergraph robustness, and found that it is related to the hypergraph degree distribution. The mulitlinear dynamical system evolution, inspired by hypergraphs, is described by the action of tensor vector multiplications between a dynamic tensor and the state vector. As a matter of fact, the multilinear dynamical system belongs to the family of homogeneous polynomial dynamical systems.

The stability of homogeneous polynomial dynamical systems is one of the most challenging problems in control theory due to its nature of nonlinearity \cite{ahmadi2019algebraic,ahmadi2013stability,ji2013constructing,samardzija1983stability,she2013discovering}. In 1983, Samardzija \cite{samardzija1983stability} established a necessary and sufficient condition for asymptotic stability in two-dimensional homogeneous polynomial dynamical systems by formulating a generalized characteristic value problem. Recently, Ali and Khadir \cite{ahmadi2019algebraic} showed that existence of a rational Lyapunov function is necessary and sufficient for asymptotic stability of a homogeneous continuously polynomial dynamical system, and the Lyapunov function can be solved using semidefinite programming. The semidefinite programming problem depends on the two degree parameters, and one has to try all possible combinations of the parameters in order to obtain the Lyapunov function. On the other hand, when a homogeneous polynomial dynamical system has degree one, the stability properties can be obtained simply from the locations of the eigenvalues of the dynamic matrix, known as the linear stability. It is therefore conceived that tensor eigenvalues may have the potential to be used for determining the stability properties of homogeneous polynomial dynamical systems. 

Tensor eigenvalue problems of real supersymmetric tensors were first explored by Qi \cite{QI20051302, QI20071363} and Lim \cite{singularvaluetensor} independently in 2005. There are many different notions of tensor eigenvalues such as H-eigenvalues, Z-eigenvalues, M-eigenvalues, and U-eigenvalues \cite{chen2021multilinear,QI20051302, QI20071363}, which are similar to matrix eigenvalues in different senses. Chen \textit{et al.} \cite{chen2017fiedler} showed that the Z-eigenvector associated with the second smallest Z-eigenvalue of a normalized Laplacian tensor can be used for hypergraph partition. In addition, Huang and Qi \cite{huang2018positive} used M-eigenvalues to prove the strong ellipticity of elasticity tensors in solid mechanics. Furthermore, Chen \textit{et al.} \cite{chen2021multilinear} utilized U-eigenvalues to determine the stability of multilinear time-invariant systems (first type of multilienar dynamical systems discussed above). Of particular interest of this paper are Z-eigenvalues. We explore the stability properties of the discrete-time version of the multilinear dynamical system proposed in \cite{chen2021controllability} via tensor spectral theory. The multilinear dynamical systems can be potential for capturing the discrete-time dynamics of hypergraphs such as coauthorship networks, film actor/actress networks, brain neural networks, and protein-protein interaction networks \cite{9119161,newman2010networks,sweeney2021network,7761624}. To the best of the author's knowledge, using tensor Z-eigenvalues to determine the stability of dynamical systems has never been reported in the literature. The key contributions of the paper are as follows:
\begin{itemize}
    \item We derive an explicit solution formula for the discrete-time multilinear dynamical systems with orthogonally decomposable (odeco) dynamic tensors by exploiting tensor Z-eigenvalues and Z-eigenvectors. 
    \item By utilizing the form of the explicit solution, we are able to discuss the stability properties of such multilinear dynamical systems. Similarly to the linear stability, the Z-eigenvalues from the orthogonal decompositions of the dynamic tensors can offer necessary and sufficient stability conditions. According to the stability conditions, we are able to obtain the regions of attraction of the multilinear dynamical systems. 
    \item We provide an upper bound for the Z-spectral radii by using tensor unfolding, which can be used to determine the asymptotic stability of the multilinear dynamical systems efficiently. 
    \item We extend the stability results to the multilinear dynamical systems with non-odeco dynamic tensors by using the tensor Frobenius norm and tensor singular values.
\end{itemize}
Note that all the results are applicable to homogeneous polynomial dynamical systems if one can find the corresponding multilinear dynamical systems. 

The paper is organized into seven sections. In Section \ref{sec:1}, we review tensor preliminaries including tensor products, tensor eigenvalues, and tensor orthogonal decomposition. In Section \ref{sec:2}, we introduce the representation of discrete-time multilinear dynamical systems. If the dynamic tensors are odeco, we can construct their explicit solutions. We establish the stability criteria for the multilinear dynamical systems with odeco dynamic tensors and propose an upper bound for the Z-spectral radii which can be used to determine the asymptotic stability efficiently in Section \ref{sec:3}. Moreover, we generalize the stability results to general multilinear dynamical systems. We verify our results with numerical examples in Section \ref{sec:4}. We discuss the stabilizability and reachability of multilinear dynamical systems with control in \ref{sec:5}, and conclude in \ref{sec:6} with future research directions.

\section{Tensor preliminaries}\label{sec:1}
A tensor is a multidimensional array \cite{9119161,doi:10.1137/1.9781611975758.18,chen2021multilinear,gelss2017tensor,Kolda06multilinearoperators, doi:10.1137/07070111X}. The order of a tensor is the number of its dimensions,  and each dimension is called a mode. A $k$th order tensor usually is denoted by $\textsf{T}\in \mathbb{R}^{n_1\times n_2\times  \dots \times n_k}$.  It is therefore reasonable to consider scalars $x\in\mathbb{R}$ as zero-order tensors, vectors $\textbf{v}\in\mathbb{R}^{n}$ as first-order tensors, and matrices $\textbf{M}\in\mathbb{R}^{m\times n}$ as second-order tensors. A tensor is called cubical if every mode is the same size, i.e., $\textsf{T}\in \mathbb{R}^{n\times n\times  \dots \times n}$. A cubical tensor \textsf{T} is called supersymmetric if $\textsf{T}_{j_1j_2\dots j_k}$ is invariant under any permutation of the indices, and is called diagonal if $\textsf{T}_{j_1j_2\dots j_k}=0$ except $j_1=j_2=\dots = j_k$.

The inner product of two tensors $\textsf{T},\textsf{S}\in \mathbb{R}^{n_1\times n_2\times \dots \times n_k}$ is defined as
\begin{equation}
\langle \textsf{T},\textsf{S}\rangle =\sum_{j_1=1}^{n_1}\dots \sum_{j_k=1}^{n_k}\textsf{T}_{j_1j_2\dots j_k}\textsf{S}_{j_1j_2\dots j_k},
\end{equation}
leading to the tensor Frobenius norm $\|\textsf{T}\|^2=\langle \textsf{T},\textsf{T}\rangle$. The tensor vector multiplication $\textsf{T} \times_{p} \textbf{v}$ along mode $p$ for a vector $\textbf{v}\in  \mathbb{R}^{n_p}$ is defined  by
\begin{equation}
(\textsf{T} \times_{p} \textbf{v})_{j_1j_2\dots j_{p-1}j_{p+1}\dots j_k}=\sum_{j_p=1}^{n_p}\textsf{T}_{j_1j_2\dots j_p\dots j_k}\textbf{v}_{j_p}, 
\end{equation}
which can be extended to 
\begin{equation}\label{eq6}
\begin{split}
\textsf{T}\times_1 \textbf{v}_1 \times_2\textbf{v}_2\times_3\dots \times_{k}\textbf{v}_k=\textsf{T}\textbf{v}_1\textbf{v}_2\dots\textbf{v}_k\in\mathbb{R}
\end{split}
\end{equation}
for $\textbf{v}_p\in \mathbb{R}^{ n_p}$. If \textsf{T} is supersymmetric and $\textbf{v}_p=\textbf{v}$ for all $p$, the product (\ref{eq6}) is also known as the homogeneous polynomial associated with \textsf{T}, and we write  it as $\textsf{T}\textbf{v}^{k}$ for simplicity. 

Homogeneous polynomials are closely related to eigenvalue problems. The tensor eigenvalues of real supersymmetric tensors were first explored by Qi \cite{QI20051302, QI20071363} and Lim \cite{singularvaluetensor} independently. There are many different notions of tensor eigenvalues such as H-eigenvalues, Z-eigenvalues, M-eigenvalues, and U-eigenvalues \cite{chen2021multilinear,QI20051302, QI20071363}. Of particular interest of this paper are Z-eigenvalues. Given a $k$th order supersymmetric tensor $\textsf{T}\in\mathbb{R}^{n\times n\times \dots \times n}$, the E-eigenvalues $\lambda\in\mathbb{C}$ and E-eigenvectors $\textbf{v}\in\mathbb{C}^n$ of \textsf{T} are defined as
\begin{equation}
\begin{cases}
   \textsf{T}\textbf{v}^{k-1} = \lambda\textbf{v}\\
   \textbf{v}^\top\textbf{v} = 1
\end{cases}.
\end{equation}
The E-eigenvalues $\lambda$ could be complex. If $\lambda$ are real, we call them Z-eigenvalues. Computing the E-eigenvalues and the Z-eigenvalues of a tensor is NP-hard \cite{Hillar:2013:MTP:2555516.2512329}. In 2016, Chen \textit{et al.} \cite{chen2016computing} proposed numerical methods for computing E-eigenvalues and Z-eigenvalues via homotopy continuation approach, but the methods only work well for small size tensors. 

There are many types of tensor decompositions including CANDECOMP/PARAFAC decomposition, higher-order singular value decomposition, Tucker decomposition, and tensor train decomposition, which all play important roles in tensor algebra \cite{doi:10.1137/S0895479896305696, Kolda06multilinearoperators, doi:10.1137/07070111X,doi:10.1137/090748330,oseledets2011tensor}.
Tensor orthogonal decomposition is a special case of  CANDECOMP/PARAFAC decomposition. A $k$th order suppersymmetric tensor $\textsf{T}\in\mathbb{R}^{n\times n\times \dots \times n}$ is called orthogonally decomposable (odeco) if it can be written as a sum of vector outer products
\begin{equation}\label{eq:22}
\textsf{T} = \sum_{r=1}^n \lambda_r \textbf{v}_r\circ\textbf{v}_r\circ \stackrel{k}{\cdots} \circ\textbf{v}_r,
\end{equation}
where $\lambda_r\in\mathbb{R}$ in the descending order, and $\textbf{v}_r\in\mathbb{R}^n$ are orthonormal \cite{doi:10.1137/140989340}. Here $``\circ"$ denotes the vector outer product. Reobeva \cite{doi:10.1137/140989340} proved that $\lambda_r$ are the Z-eigenvalues of \textsf{T} with the corresponding Z-eigenvectors $\textbf{v}_r$. However, $\lambda_r$ do not include all the Z-eigenvalues of \textsf{T}. Moreover, the author showed that odeco tensors satisfy a set of polynomial equations that vanish on the odeco variety, which is the Zariski closure of the set of odeco tensors inside the space of $k$th order $n$-dimensional complex supersymmetric tensors. Note that the author only proved for the case when $n=2$, but provided with strong evidence for its overall correctness \cite{doi:10.1137/140989340}. A tensor power method was also reported in \cite{doi:10.1137/140989340} in order to find the orthogonal decomposition of an odeco tensor.

\section{Multilinear Dynamical Systems}\label{sec:2}
In this paper, we consider the discrete-time version of the multilinear dynamical system proposed in \cite{chen2021controllability}, which is given by
\begin{equation}\label{eq:1}
\textbf{x}_{t+1} = \textsf{A}\times_1\textbf{x}_t\times_2\textbf{x}_t\times_3\dots\times_{k-1}\textbf{x}_t= \textsf{A}\textbf{x}_t^{k-1},
\end{equation}
where $\textsf{A}\in\mathbb{R}^{n\times n\times \dots \times n}$ is a supersymmetric dynamic tensor, and $\textbf{x}_t\in\mathbb{R}^n$ is the state variable (multilinear in the sense of multilinear algebra). In fact, the multilinear dynamical system (\ref{eq:1}) belongs to the family of homogeneous polynomial dynamical systems of degree $k-1$ if one expands the tensor vector multiplications. 

\subsection{Explicit Solutions}
Finding an explicit solution of a multilinear dynamical system is usually challenging due to its nonlinear nature. However, if the dynamic tensor \textsf{A} is odeco, we can write down the solution of (\ref{eq:1}) explicitly in a simple fashion. 

\textit{Proposition 1:} Suppose that $k\geq 3$ and $\textsf{A}\in\mathbb{R}^{n\times n\times \dots \times n}$ is odeco with orthogonal decomposition (\ref{eq:22}). Let the initial condition $\textbf{x}_0=\sum_{r=1}^nc_r\textbf{v}_r$. Then the solution of the multilinear dynamical system (\ref{eq:1}) at time $q$, given initial condition $\textbf{x}_0$, is given by
\begin{equation}
    \textbf{x}_q = \sum_{r=1}^n\lambda_r^{\alpha}c_r^{\beta}\textbf{v}_r,
\end{equation}
where $\alpha=\sum_{j=0}^{q-1}(k-1)^j=\frac{(k-1)^q-1}{k-2}$ and $\beta = (k-1)^q$.

\textit{Proof:} Based on the properties of tensor vector multiplications and tensor orthongonal decomposition, we can write down the solution $\textbf{x}_1$ as follows:
\begin{align*}
    \textbf{x}_1 &= \textsf{A}\times_1\Big(\sum_{i=1}^nc_i\textbf{v}_i\Big)\times_2 \Big(\sum_{i=1}^nc_i\textbf{v}_i\Big)\times_3\dots\times_{k-1}\Big(\sum_{i=1}^nc_i\textbf{v}_i\Big)\\
    & = \Big(\sum_{r=1}^n \lambda_r \textbf{v}_r\circ\textbf{v}_r\circ \dots \circ\textbf{v}_r\Big)\times_1\Big(\sum_{i=1}^nc_i\textbf{v}_i\Big)\times_2 \Big(\sum_{i=1}^nc_i\textbf{v}_i\Big)\\ &\times_3  \dots\times_{k-1}\Big(\sum_{i=1}^nc_i\textbf{v}_i\Big)= \sum_{r=1}^n \lambda_r \Big\langle \textbf{v}_r, \sum_{i=1}^nc_i\textbf{v}_i\Big\rangle^{k-1} \textbf{v}_r.
\end{align*}
Since all the vectors $\textbf{v}_r$ are orthonormal, we have
\begin{equation*}
    \textbf{x}_1=\sum_{r=1}^n \lambda_r \Big\langle \textbf{v}_r, \sum_{i=1}^nc_i\textbf{v}_i\Big\rangle^{k-1} \textbf{v}_r=\sum_{r=1}^n \lambda_r\textbf{c}_r^{k-1}\textbf{v}_r.
\end{equation*}
Similarly, the solution $\textbf{x}_2$ can be written as 
\begin{align*}
    \textbf{x}_2 &= \textsf{A}\times_1\Big(\sum_{i=1}^n\lambda_ic_i^{k-1}\textbf{v}_i\Big)\times_2 \Big(\sum_{i=1}^n\lambda_ic_i^{k-1}\textbf{v}_i\Big)\\ &\times_3\dots\times_{k-1}\Big(\sum_{i=1}^n\lambda_ic_i^{k-1}\textbf{v}_i\Big)\\
    & = \sum_{r=1}^n \lambda_r \Big\langle \textbf{v}_r, \sum_{i=1}^n\lambda_ic_i^{k-1}\textbf{v}_i\Big\rangle^{k-1} \textbf{v}_r = \sum_{r=1}^n \lambda_r^k\textbf{c}_r^{(k-1)^2}\textbf{v}_r.
\end{align*}
One can continue to compute $\textbf{x}_3,\textbf{x}_4,\dots,\textbf{x}_q$ in the similar manner. Therefore, the result follows immediately. \hfill  $\blacksquare$

The coefficient $c_r$ can be found from the inner product between $\textbf{x}_0$ and $\textbf{v}_r$, and thus one may write $|c_r\lambda_r^{\frac{1}{k-2}}|$ as $|\langle \textbf{x}_0,\lambda_r^{\frac{1}{k-2}}\textbf{v}_r\rangle|$. When $k=2$, Proposition 1 reduces to the explicit solution formula for linear dynamical systems, i.e., 
\begin{equation*}
    \lim_{k\rightarrow 2} \lambda_r^{\frac{(k-1)^q-1}{k-2}}c_r^{(k-1)^q}=\lim_{k\rightarrow 2} \lambda_r^{q(k-1)^{q-1}}c_r^{(k-1)^q}=c_r\lambda_r^q.
\end{equation*}
Moreover, by exploiting the form of the explicit solution, we are able to establish the stability criteria for the multilinear dynamical systems with odeco dynamic tensors. 

\section{Stability} \label{sec:3}
In linear control theory, it is conventional to investigate so-called (internal) stability. The stability of a linear dynamical system solely depends on the locations of the eigenvalues of the dynamic matrix. Similarly to linear dynamical systems, the equilibrium point $\textbf{x}=\textbf{0}$ of the multilinear dynamical system (\ref{eq:1}) is called stable if $\|\textbf{x}_t\|\leq \gamma \|\textbf{x}_0\|$ for some initial condition $\textbf{x}_0$ and $\gamma>0$, asymptotically stable if $\textbf{x}_t\rightarrow \textbf{0}$ as $t\rightarrow \infty$, and unstable if it is not stable. 

The stability properties of the multilinear dynamical system (\ref{eq:1}) with an odeco dynamic tensor are similar to these of linear dynamical systems, but the stability depends on both the Z-eigenvalues of the dynamic tensor \textsf{A} and initial conditions.

\textit{Proposition 2}
Suppose that $k\geq 3$ and $\textsf{A}\in\mathbb{R}^{n\times n\times \dots \times n}$ is odeco with orthogonal decomposition (\ref{eq:22}). Let the initial condition $\textbf{x}_0=\sum_{r=1}^nc_r\textbf{v}_r$. For the multilinear dynamical system (\ref{eq:1}), the equilibrium point $\textbf{x}=\textbf{0}$ is: 
\begin{itemize}
    \item stable if and only if $|c_r\lambda_r^{\frac{1}{k-2}}|\leq 1$ for all $r=1,2,\dots,n$;
    \item asymptotically stable if and only if $|c_r\lambda_r^{\frac{1}{k-2}}|< 1$ for all $r=1,2,\dots,n$;
    \item unstable if and only if $|c_r\lambda_r^{\frac{1}{k-2}}|> 1$ for some $r=1,2,\dots,n$.
\end{itemize}

\textit{Proof:} Based on the result from Proposition 1, the solution at time $q$, given initial condition $\textbf{x}_0$, is $\textbf{x}_q = \sum_{r=1}^n\lambda_r^{\alpha}c_r^{\beta}\textbf{v}_r$ where $\alpha=\sum_{j=0}^{q-1}(k-1)^j=\frac{(k-1)^q-1}{k-2}$ and $\beta = (k-1)^q$. Therefore, it can be shown that
\begin{equation*}
    \lambda_r^{\alpha}c_r^{\beta} = \lambda_r^{\frac{(k-1)^q-1}{k-2}}c_r^{(k-1)^q}=\lambda_r^{-\frac{1}{k-2}}(\lambda_r^{\frac{1}{k-2}}c_r)^{(k-1)^q}.
\end{equation*}
Hence, the results follow immediately.  \hfill $\blacksquare$

The inequalities obtained from the asymptotic stability condition  can provide us with the region of attraction of the multilinear dynamical system (\ref{eq:1}), i.e., 
\begin{equation}
    R = \{\textbf{x}:|c_r|<|\lambda_r|^{-\frac{1}{k-2}} \text{ where }\textbf{x}=\sum_{r=1}^nc_r\textbf{v}_r\}.
\end{equation}
Moreover, if the product between $\max_r{|c_r|}$ and $\max_r{|\lambda_r|^{\frac{1}{k-2}}}$ is less than one, the multilinear dynamical system (\ref{eq:1}) will be asymptotically stable. 

\textit{Definition 1:} The Z-spectral radius of a supersymmetric tensor is the maximum of the absolute values of all its Z-eigenvalues.

\textit{Corollary 1:} Suppose that $k\geq 3$ and $\textsf{A}\in\mathbb{R}^{n\times n\times \dots \times n}$ is odeco. Let $\textbf{x}_0$ be some initial conditions. For the multilinear dynamical system (\ref{eq:1}), the equilibrium point $\textbf{x}=\textbf{0}$ is asymptotically stable if $\lambda^{\frac{1}{k-2}}\|\textbf{x}_0\|<1$ where $\lambda$ is the Z-spectral radius of \textsf{A}.  

\textit{Proof:} By the Cauchy-Schwarz inequality, $|c_r|\leq \|\textbf{x}_0\|$ for all $r$. In addition, $\max_r{|\lambda_r|}\leq \lambda$. Therefore, the result follows immediately from Proposition 2. \hfill $\blacksquare$

Based on our numerical experiments, we find that the Z-spectral radius 
\begin{equation*}
    \lambda=\max{\{|\lambda_1|,|\lambda_n|\}}, 
\end{equation*}
where $\lambda_1$ and $\lambda_n$ are the largest and the smallest Z-eigenvalues in the orthogonal decomposition, respectively. This implies that $\lambda_1$ is the largest Z-eigenvalue of \textsf{A}, or $\lambda_n$ is the smallest Z-eigenvalue of \textsf{A}. However, the correctness of this conjecture needs further investigation. 

\subsection{Upper Bounds for Z-spectral Radii}
Computing the orthogonal decomposition or Z-eigenvalues of a supersymmetric tensor is NP-hard even if we know the tensor is odeco beforehand (e.g., the tensor satisfies a set of polynomial equations that vanish on the odeco variety) \cite{Hillar:2013:MTP:2555516.2512329,doi:10.1137/140989340}. If we can come up with some upper bounds for the Z-spectral radii of the dynamic tensors, it will save a great amount of computations to determine the stability of the multilinear dynamical systems with odeco dynamic tensors.

\textit{Lemma 1:} Suppose that $\textsf{A}\in\mathbb{R}^{n\times n\times \dots \times n\times n}$ is an even-order supersymmetric tensor. The Z-spectral radius of \textsf{A} is upper bounded by the spectral radius of the unfolded matrix defined by the following
\begin{equation}\label{eq:3}
    \textbf{A} = \psi(\textsf{A}) \text{ s.t. } \textsf{A}_{j_1i_1\dots j_ki_k} \xrightarrow{\psi} \textbf{A}_{ji},
\end{equation}
where $j=j_1+\sum_{p=2}^k(j_p-1)n^{p-1}$, and $i=i_1+\sum_{p=2}^k(i_p-1)n^{p-1}$. 

The essence of the proof is based on the fact that the restriction of $\psi^{-1}$ on the general linear group produces a group isomorphism \cite{doi:10.1137/1.9781611975758.18,chen2021multilinear}. We omit detailed proof in this paper. The eigenvalues of the unfolded matrix \textbf{A} are also called the U-eigenvalues of \textsf{A} \cite{doi:10.1137/1.9781611975758.18,chen2021multilinear}. The computational cost for computing spectral radii from the matrix eigenvalue decomposition is much less than that for Z-spectral radii. Once we have an upper bound for the Z-spectral radii of the dynamic tensors, we can determine the asymptotic stability of the multilinear dynamical system (\ref{eq:1}) without computing orthogonal decomposition or Z-eigenvalues.

\textit{Corollary 2:} Suppose that $k\geq 4$ is even and $\textsf{A}\in\mathbb{R}^{n\times n\times \dots \times n\times n}$ is odeco. Let $\textbf{x}_0$ be some initial conditions. For the multilinear dynamical system (\ref{eq:1}), the equilibrium point $\textbf{x}=\textbf{0}$ is asymptotically stable if $\mu^{\frac{1}{k-2}}\|\textbf{x}_0\|<1$ where $\mu$ is the spectral radius of $\psi(\textsf{A})$.

\textit{Proof:} The result follows immediately from Lemma 1 and Corollary 1.  \hfill    $\blacksquare$

The condition offers a conservative region of attraction for the multilinear dynamical system (\ref{eq:1}) without knowing the orthogonal decomposition of \textsf{A}, i.e., 
\begin{equation*}
R = \{\textbf{x}:\|\textbf{x}\|< \mu^{-\frac{1}{k-2}}\}.
\end{equation*}
In addition, there are many other upper bounds for Z-spectral radii of supersymmetric tensors \cite{chang2013some,he2014upper,ma2019some,wu2018upper}. For example, He \textit{et al.} \cite{he2014upper} proposed that given a positive $k$th order supersymmetric tensor \textsf{A}, its Z-spectral radius is upper bounded by 
\begin{equation}\label{eq:8}
    \lambda \leq R - l\Big(1-(\frac{r}{R})^{\frac{1}{k}}\Big),
\end{equation}
where $l$ is the minimum entry of \textsf{A},  
\begin{align*}
r&=\min_{j}\Big(\sum_{j_2=1}^n\dots \sum_{j_k=1}^n\textsf{A}_{jj_2\dots j_k}\Big),\\
R&=\max_{j}\Big(\sum_{j_2=1}^n\dots \sum_{j_k=1}^n\textsf{A}_{jj_2\dots j_k}\Big).
\end{align*}
Hence, one can also use this upper bound to determine the stability of a multilinear dynamical system if the dynamic tensor contains all positive entries. Note that given a dynamic tensor, the better upper bound of the Z-spectral radius, the more strong stability conditions we can obtain.

\subsection{Non-odeco Dynamic Tensors} \label{sec:4}
As mentioned in \cite{doi:10.1137/140989340}, not all supersymmetric tensors are odeco. Therefore, we offer a more general but relatively weaker stability result for the multilinear dynamical system (\ref{eq:1}) with non-odeco dynamic tensors. First, we adapt a lemma from \cite{doi:10.1002/nla.2086} that provides a bound on the Frobenius norm of tensor vector multiplications.

\textit{Lemma 2:} Suppose that $\textsf{A}\in\mathbb{R}^{n\times n\times \dots \times n}$ and $\textbf{v}\in\mathbb{R}^n$. Then 
\begin{equation*}
    \|\textsf{A}\times_p\textbf{v}\|\leq \|\textsf{A}\|\|\textbf{v}\|.
\end{equation*}

\textit{Proposition 3:} Suppose that $k\geq 3$. Let $\textbf{x}_0$ be some initial conditions. For the multilinear dynamical system (\ref{eq:1}), the equilibrium point $\textbf{x}=\textbf{0}$ is asymptotically stable if 
\begin{equation*}
  \|\textsf{A}\|^{\frac{1}{k-2}}\|\textbf{x}_0\|< 1.
\end{equation*}

\textit{Proof:} Based on Lemma 2, we have 
\begin{equation*}
    \|\textbf{x}_{t+1}\|\leq \|\textsf{A}\|\|\textbf{x}_t\|^{k-1}.
\end{equation*}
Thus, it can be shown similarly as Proposition 1 that at the $q$th step, we have
\begin{equation*}
    \|\textbf{x}_q\|\leq \|\textsf{A}\|^\alpha\|\textbf{x}_0\|^\beta, 
\end{equation*}
where $\alpha$ and $\beta$ are the same quantities as defined in Proposition 1. Therefore, the result follows immediately. \hfill $\blacksquare$

Proposition 3 also holds for non-supersymmetric dynamic tensors \textsf{A}. Similarly, we can obtain a conservative region of attraction, i.e., 
\begin{equation*}
    R = \{\textbf{x}:\|\textbf{x}\|< \|\textsf{A}\|^{-\frac{1}{k-2}}\}.
\end{equation*}
Moreover, the tensor Frobenius norm is closely related to generalized $p$-mode singular values. Detailed definitions of generalized $p$-mode singular values can be found in \cite{doi:10.1137/S0895479896305696}.

\textit{Corollary 3:} Suppose that $k\geq 3$. Let $\textbf{x}_0$ be some initial conditions. For the multilinear dynamical system (\ref{eq:1}), the equilibrium point $\textbf{x}=\textbf{0}$ is asymptotically stable if for any $p$,
\begin{equation*}
  \Big(\sum_{j=1}^n(\gamma_j^{(p)})^2\Big)^{\frac{1}{k-2}}\|\textbf{x}_0\|< 1, 
\end{equation*}
where $\gamma_j^{(p)}$ are the $p$-mode singular values of \textsf{A}.

\textit{Proof:} The result follows immediately from Proposition 3 and the fact that the Frobenius norm of a tensor is equal to the sum of its $p$-mode singular values' square for any $p$. \hfill $\blacksquare$

\section{Numerical Examples:} \label{sec:4}
All the numerical examples presented were performed on a Macintosh machine with 16 GB RAM and a 2 GHz Quad-Core Intel Core i5 processor in MATLAB R2020b, and used the MATLAB tensor toolbox \cite{BaKo06,tensortoolbox} and MATLAB TT toolbox \cite{oseledets2012tt,oseledets2011tensor}.

\begin{table}[t]
\caption{Stability results for different initial conditions.}
\centering
\begin{tabular}{|c|c|c|c|}
\hline
\textbf{IC}  & $\max{\{|c_r\lambda_r|\}}$ & $\|\textsf{A}\|\|\textbf{x}_0\|$&\textbf{Stability} \\ \hline
 \textbf{a} &  0.9735  & 28.7712 & Asym. Stable\\ \hline
 \textbf{b}     &   0.6032     &   0.9413       &    Asym. Stable \\ \hline
 \textbf{c}     &    1     &    53.9410     &  Stable                                   \\ \hline
 \textbf{d}    &       1.0053      &  1.5688   &      Unstable   \\ \hline
\end{tabular}
\label{tab:3}
\end{table}

\definecolor{bleudefrance}{rgb}{0.19, 0.55, 0.91}
\definecolor{emerald}{rgb}{0.31, 0.78, 0.47}
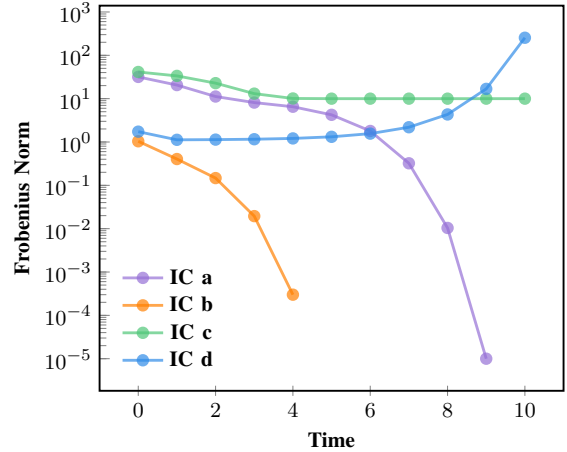
\begin{figure}[t!]
\centering
\begin{tikzpicture}[scale=0.9]
\begin{axis}[ymode=log, ytick={1e-5, 1e-4, 1e-3, 1e-2, 1e-1, 1, 10, 100, 1000}, xlabel={\textbf{Time}}, ylabel={\textbf{Frobenius Norm}}, ylabel near ticks, legend pos=south west, xlabel near ticks,xtick pos=left,ytick pos=left,xtick align=center,legend style={draw=none},axis line style = thick]
\addplot+[very thick,mark=*,opacity=.7,mypurple,mark options={draw=mypurple,fill=mypurple}] plot coordinates
{ (0,31.7648) (1,20.5942) (2,11.2037) (3, 8.1222) (4,6.5110) (5,4.2388) (6,1.7968) (7,0.3228) (8,0.0104) (9,0.00001) (10,0)};
\addlegendentry{\textbf{IC a}}
\addplot+[very thick,mark=*,opacity=.7, orange,mark options={draw=orange,fill=orange}] plot coordinates
{ (0,1.0392) (1,0.4045) (2,0.1471) (3, 0.0195) (4,0.0003) (5,0) (6,0) (7,0) (8,0) (9,0) (10,0)};
\addlegendentry{\textbf{IC b}}
\addplot+[very thick,mark=*,opacity=.7,emerald,mark options={draw=emerald,fill=emerald}] plot coordinates
{ (0,41.2432) (1,33.5382) (2,22.8027) (3, 13.0613) (4,10.1006) (5,10.0002) (6,10) (7,10) (8,10) (9,10) (10,10) };
\addlegendentry{\textbf{IC c}}
\addplot+[very thick,mark=*,opacity=.7, bleudefrance,mark options={draw=bleudefrance,fill=bleudefrance}] plot coordinates
{ (0,1.7321) (1,1.1235) (2,1.1349) (3,1.1593) (4,1.2096) (5,1.3167) (6,1.5604) (7,2.1914) (8,4.3221) (9,16.8127) (10,254.4007)};
\addlegendentry{\textbf{IC d}}
\end{axis}
\end{tikzpicture}
\hspace{0.9cm}
\caption{Stability results for different initial conditions corresponding to Table \ref{tab:3}. When the norm of $\textbf{x}_t$ is less than $10^{-5}$, we omit the point.}
\label{fig:12}
\end{figure}

\subsection{Stability with Odeco Dynamic Tensors}
In this example, we try to verify the stability results discussed in Proposition 2. Given an odeco tensor $\textsf{A}\in\mathbb{R}^{3\times 3\times 3}$ with the orthogonal decomposition (columns of \textbf{V} are $\textbf{v}_r$ in (\ref{eq:22}))
\begin{equation*}
    \textbf{V} = \begin{bmatrix} -0.8482  & -0.5212  &  0.0947\\-0.4840  &  0.6899 &  -0.5382\\ 0.2152 &  -0.5024 &  -0.8374\end{bmatrix} \text{ and } \boldsymbol{\lambda} = \begin{bmatrix} 0.9\\0.1\\0.02\end{bmatrix},
\end{equation*}
we compute the trajectories $\textbf{x}_t$ for four initial conditions, which are given by
\begin{equation*}
    \textbf{x}_\textbf{a} = \begin{bmatrix} 3 \\10\\30\end{bmatrix}, \text{ } \textbf{x}_\textbf{b} = \begin{bmatrix} 0.6 \\0.6\\0.6\end{bmatrix}, \text{ } \textbf{x}_\textbf{c} = \begin{bmatrix} -2.2720\\-15.1148\\-38.3064\end{bmatrix}, \text{ } \textbf{x}_\textbf{d} = \begin{bmatrix} 1 \\1\\1\end{bmatrix}.
\end{equation*}
The results are shown in Table \ref{tab:3} and  Figure \ref{fig:12}. For each initial condition, we calculate the quantities $\max{\{|c_r\lambda_r|\}}$ and $\|\textsf{A}\|\|\textbf{x}_0\|$, and compare them to one. Clearly, the locations of $c_r\lambda_r$ determine the stability of the multilinear dynamical system. The  region of attraction $R$ of the multilinear dynamical system can be obtained by
\begin{equation*}
    R = \left\{\textbf{x}:\begin{array}{lr}
        |-0.8482x_1-0.4840x_2+0.2152x_3|<\frac{10}{9}\\
        |-0.5212x_1+0.6899x_2-0.5024x_3|<10\\
        |0.0947x_1-0.5382x_2-0.8374x_3|<50
        \end{array}\right\},
\end{equation*}
where $\textbf{x}=\begin{bmatrix}x_1 & x_2 & x_3\end{bmatrix}^\top$. In addition, the stability condition stated in Proposition 3 is weaker than that in Proposition 2 as seen in Figure \ref{fig:12} IC a and b.

\subsection{Stability Using Upper Bounds of Z-spectral Radii}
\begin{figure}[t]
\centering
\begin{tikzpicture}[scale=1.3]
\begin{axis}[ymin=-2.2,ymax=2.2,xmax=2.2,xmin=-2.2,xtick=\empty, ytick=\empty, axis equal image,axis lines = middle,axis line style = thick,legend pos=south west,legend style={draw=none}]
\addplot+[very thick,mark=*,opacity=.7,bleudefrance,mark options={draw=bleudefrance,fill=bleudefrance}] plot coordinates
{ (1.2,1.2) (1.7737,1.6769) (5.2838,4.9951)};

\addplot+[very thick,mark=*,opacity=.7,orange,mark options={draw=orange,fill=orange}] plot coordinates
{ (1,1) (1.0264,0.9704) (1.0240,0.9681) (1.0168,0.9612) (0.9954) (0.9410) (0.9338,0.8828) (0.7710,0.7288) (0.4339,0.4102) (0.0773,0.0731) (0.0004381,0.0004142)};

\addplot+[very thick,mark=*,opacity=.7,mypurple,mark options={draw=mypurple,fill=mypurple}] plot coordinates
{ (-1.4,0) (-0.6270, -0.1031) (-0.0654,-0.0370) (-0.0001428,-0.0001317)};

\addplot+[very thick,mark=*,opacity=.7,emerald,mark options={draw=emerald,fill=emerald}] plot coordinates
{ (0.9,-0.9) (0.5659, -0.5986) (0.1536,-0.1625) (0.0031,-0.0033) (0,0)};

\node [below] at (axis cs: 2.1, 0) {$x_1$};
\node [right] at (axis cs: 0, 2.1) {$x_2$};
\node [below] at (axis cs: 1.6, 0) {$\sqrt{2}$};
\node [right] at (axis cs: 0, 1.6) {$\sqrt{2}$};
\node [below] at (axis cs: -1.85, 0) {$-\sqrt{2}$};
\node [right] at (axis cs: 0, -1.6) {$-\sqrt{2}$};
\node [right] at (axis cs: 1.3, 1.2) {\small{\textbf{IC d}}};
\node [above] at (axis cs: -1.8, 0) {\small{\textbf{IC a}}};
\node [left] at (axis cs: 0.85, -0.9) {\small{\textbf{IC b}}};
\node [left] at (axis cs: 0.9, 1) {\small{\textbf{IC c}}};
\draw [style=ultra thick, red, radius=1.4142, dashed] (axis cs:0,0) circle;
\end{axis}
\end{tikzpicture}
\hspace{0.9cm}
\caption{A conservative region of attraction (red dashed line) of the multilinear dynamical system, and four initial conditions with their trajectories.}
\label{fig:2}
\end{figure}
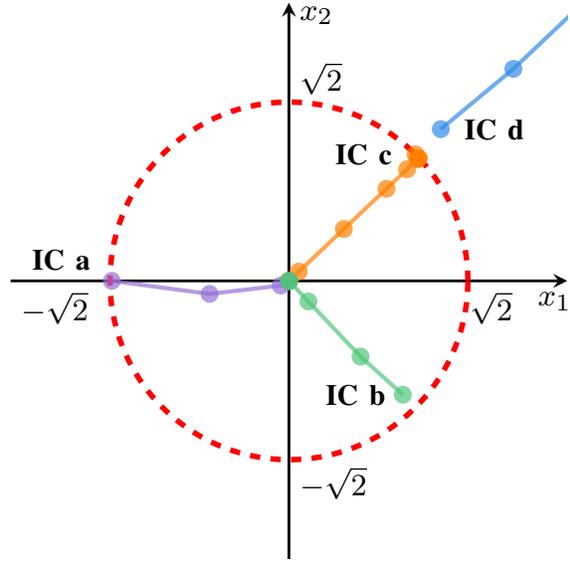

In this example, we try to apply the upper bound of the Z-spectral radii defined in (\ref{eq:3}) to obtain a conservative region of attraction for a multilinear dynamical system. Suppose that the dynamic tensor $\textsf{A}\in\mathbb{R}^{2\times 2\times 2\times 2}$ is odeco and is given by
\begin{align*}
    &\small{\textsf{A}_{::11} = \begin{bmatrix}0.2285 & 0.0376\\ 0.0376 & 0.2243\end{bmatrix}}, \text{ } \small{\textsf{A}_{::12} = \begin{bmatrix}0.0376 & 0.2243\\ 0.2243 & 0.0124\end{bmatrix}},\\
    &\small{\textsf{A}_{::21} = \begin{bmatrix}0.0376 & 0.2243\\ 0.2243 & 0.0124\end{bmatrix}}, \text{ } \small{\textsf{A}_{::22} = \begin{bmatrix}0.2243  & 0.0124\\ 0.0124 & 0.2229\end{bmatrix}}.
\end{align*}
The unfolded matrix $\psi(\textsf{A})$ therefore is given by
\begin{equation*}
    \psi(\textsf{A}) = \begin{bmatrix}0.2285  &  0.0376 &   0.0376  &  0.2243\\0.0376  &  0.2243 &   0.2243  &  0.0124\\0.0376 &   0.2243  &  0.2243  &  0.0124\\0.2243 &   0.0124 &   0.0124   & 0.2229\end{bmatrix}.
\end{equation*}
The spectral radius of $\psi(\textsf{A})$ is $\mu=\frac{1}{2}$, and thus the conservative region of attraction of the multilinear dynamical system is an open disk with radius $\sqrt{2}$ centered at the origin (note that the second upper bound (\ref{eq:8}) produces $\lambda\leq 1.0263$, which will give an even more conservative region of attraction). We test four initial conditions to verify the region of attraction, which are given by 
\begin{equation*}
    \textbf{x}_\textbf{a} = \begin{bmatrix}-1.4 \\0\end{bmatrix}, \text{ } \textbf{x}_\textbf{b} = \begin{bmatrix} 0.9 \\-0.9\end{bmatrix}, \text{ } \textbf{x}_\textbf{c} = \begin{bmatrix} 1\\1\end{bmatrix}, \text{ } \textbf{x}_\textbf{d} = \begin{bmatrix} 1.2\\1.2\end{bmatrix}.
\end{equation*}
The results are shown in Figure \ref{fig:2}. It is clear to see that the trajectories of the multilinear dynamical system with initial conditions started within the open disk converge to the origin, see IC a and b. Moreover, since the region of attraction is conservative, we see a trajectory started on the circle also converge to the origin, see IC c.

\section{Discussion}\label{sec:5}
The  two numerical examples reported here highlight that tensor spectral theory plays a significant role in the stability analysis of the discrete-time multilinear dynamical system (\ref{eq:1}). In particular, when the dynamic tensor \textsf{A} is odeco, its Z-eigenvalues together with initial conditions can provide necessary and sufficient criteria to determine the stability of a multilinear dynamical system. However, more theoretical and numerical investigations are required to evaluate the stability properties of multilinear dynamical systems with non-odeco dynamic tensors. Moreover, stabilizability of the multilinear dynamical system (\ref{eq:1}) needs to be considered for control/feedback designs. Stabilizability is also associated to another important concept called reachability in the systems theory. 

We are interested in exploring the reachability of the multilinear dynamical system (\ref{eq:1}) with linear control inputs, i.e., 
\begin{equation}\label{eq:7.1}
    \textbf{x}_{t+1} = \textsf{A}\textbf{x}_t^{k-1}+\textbf{B}\textbf{u}_t,
\end{equation}
where $\textbf{B}\in\mathbb{R}^{n\times m}$ is the control matrix, and $\textbf{u}_t\in\mathbb{R}^{m}$ is a control input. Chen \textit{et al.} \cite{chen2021controllability} proved that for the continuous-time case, i.e., the multilinear dynamical system is controllable if and only if $k$ is even, and the span of the smallest Lie algebra of vector fields $\{\textsf{A},\textbf{b}_1,\textbf{b}_2,\dots,\textbf{b}_m\}$ is $\mathbb{R}^n$ at the origin where $\textbf{B}=\begin{bmatrix}\textbf{b}_1 & \textbf{b}_2 & \dots & \textbf{b}_m\end{bmatrix}$. We believe that a similar result should hold for the discrete-time case, analogous to that in linear systems theory.   

\textit{Definition 2:} Suppose that $\textbf{B}=\begin{bmatrix}\textbf{b}_1 & \textbf{b}_2 & \dots & \textbf{b}_m\end{bmatrix}$. Let $\mathscr{R}_0$ be the linear span of the vectors  $\{\textbf{b}_1,\textbf{b}_2,\dots,\textbf{b}_m\}$ and $\textsf{A}\in\mathbb{R}^{n\times n\times \dots\times n}$ be a supersymmetric tensor. For each integer $q\geq 1$, define $\mathscr{R}_q$ inductively as the linear span of 
\begin{equation}
\mathscr{R}_{q-1}\cup \{\textsf{A}\textbf{v}_1\textbf{v}_2\dots\textbf{v}_{k-1}|\textbf{v}_l\in \mathscr{R}_{q-1}\}.
\end{equation}
Denote the subspace $\mathscr{R}=\cup_{q\geq 0} \mathscr{R}_q$. 

\textit{Conjecture 1:} Suppose that $k$ is even. The multilinear dynamical system (\ref{eq:7.1}) is reachable if and only if the subspace $\mathscr{R}$ spans $\mathbb{R}^n$, or equivalently, the matrix \textbf{R}, including all the column vectors from $\mathscr{R}$, has rank $n$. 

Although we are not able to fully prove Conjecture 1, one can validate its correctness with help of Gr$\ddot{\text{o}}$bner basis for specific values of $n$ and $k$. In addition, the reason that $k$ has to be even is because the roots of polynomial systems of even degree might all be complex. 

\section{Conclusion} \label{sec:6}

In this paper, we investigated the stability properties of discrete-time multilinear dynamical systems, inspired by hypergraphs. In contrast to  linear dynamical systems, the stability of the multilinear dynamical system (\ref{eq:1}) depends on both the spectrum of the dynamic tensor \textsf{A} and initial conditions. In particular, when the dynamic tensor $\textsf{A}$ is odeco, we can obtain necessary and sufficient conditions by exploiting tensor Z-eigenvalues. We also provided an upper bound for the Z-spectral radii of even-order supersymmetric tensors, which can be used to determine the asymptotic stability of the multilinear dynamical system efficiently. In addition, we extend the stability results to the multilinear dynamical systems with non-odeco dynamic tensors by using the tensor Frobenius norm and $p$-mode singular values. All the results are applicable to homogeneous polynomial dynamical systems if one can find the corresponding multilinear dynamical systems. 

It will be worthwhile to explore more strong stability conditions regarding the multilinear dynamical systems with non-odeco dynamic tensors and to realize the Lyapunov theory in multilinear dynamical systems. For example, how does one derive the Lyapunov equations for the multilinear dynamical systems (\ref{eq:1})? It will also be interesting to explore the stability properties of continuous-time multilinear dynamical systems. As mentioned in Section \ref{sec:5}, we also intend to analyze the stabilizability and reachability of multilinear control systems in future work. One particular application we plan to investigate is that of higher-order genomic networks. Recent advances in genomics technology, such as multiway chromosomal conformation capture (Pore-C), have inspired us to consider the human genome as a hypergraph \cite{ulahannan2019nanopore}. Stability and stabilizability will be important in analyzing such higher-order networked systems.

\section*{Acknowledgements}
I would like to thank my advisor Dr. Anthony M. Bloch for carefully reading the manuscript and for providing constructive comments.

\bibliographystyle{IEEEtran}
\bibliography{reference}

\end{document}